\documentclass[11pt]{amsproc}
\usepackage[mathcal]{euscript}
\usepackage{amsthm}
\usepackage{amsfonts}
\usepackage{amsmath,amssymb}
\usepackage{indentfirst,latexsym,bm,amsthm,graphicx,colortbl}
\usepackage{fancyhdr}
\usepackage{times}
\usepackage{amsmath,amsfonts,amssymb,amsthm}
\usepackage{epsfig}
\usepackage[all]{xy}
\usepackage{tikz}
\usepackage{color}
\usepackage{enumerate}
\usepackage{verbatim}
\newtheorem{theorem}{Theorem}[section]

\newtheorem{prop}[theorem]{Proposition}

\theoremstyle{definition}
\newtheorem{defi}[theorem]{Definition}
\newtheorem{example}[theorem]{Example}

\newtheorem{coro}[theorem]{Corollary}
\theoremstyle{remark}
\newtheorem{remark}[theorem]{Remark}

\numberwithin{equation}{section}

\def \<{\langle}
\def \>{\rangle}

\def \be{\begin{equation}\label}
\def \ee{\end{equation}}
\def \bex{\begin{example}\label}
\def \eex{\end{example}}
\def \bl{\begin{lem}\label}
\def \el{\end{lem}}
\def \bt{\begin{thm}\label}
\def \et{\end{thm}}
\def \bp{\begin{prop}\label}
\def \ep{\end{prop}}
\def \br{\begin{rem}\label}
\def \er{\end{rem}}
\def \bc{\begin{coro}\label}
\def \ec{\end{coro}}
\def \bd{\begin{de}\label}
\def \ed{\end{de}}

\begin{document}
\title[Vertex representations of quantum $N$-toroidal algebras for type $C$ ]
{Vertex representations of quantum $N$-toroidal algebras for type $C$}


\author[Jing]{Naihuan Jing}
\address{Department of Mathematics, Shanghai University, Shanghai 200444, China\newline \indent
Department of Mathematics, North Carolina State University,
   Raleigh, NC 27695, USA \\}
\email{jing@math.ncsu.edu}

\author[Xu]{Zhucheng Xu}
\address{Department of Mathematics, Shanghai University,
Shanghai 200444, China} 

\author[Zhang]{Honglian Zhang$^\star$}
\address{Department of Mathematics, Shanghai University,
Shanghai 200444, China} \email{hlzhangmath@shu.edu.cn}


\begin{abstract}
Quantum $N$-toroidal algebras are generalizations of quantum affine algebras and quantum toroidal algebras. In this paper 
we construct a level-one vertex representation of the quantum $N$-toroidal algebra for type $C$. In particular, we also obtain a
level-one module of the quantum toroidal algebra for type $C$ as a special case.
\end{abstract}

\maketitle

\section{Introduction}

The $N$-toroidal Lie algebra is the universal central extension of the loop algebra $\mathfrak{g}\otimes \mathbb{C}[t_1^{\pm1},$
$ \cdots, t_{N}^{\pm1}]$ associated to the simple Lie algebra $\mathfrak{g}$. The most important subclasses are the 1-toroidal or the affine Kac-Moody algebra $\widehat{\mathfrak g}$ and the 2-toroidal Lie algebra or simply the toroidal Lie algebra.
It is well-known that the affine Kac-Moody algebra is the one-dimensional universal central extension of the loop algebra, while the $N$-toroidal Lie algebra ($N\geq 2)$ is an infinite-dimensional central extension of the loop algebra. Extensive studies have been done on
representations of toroidal Lie algebras (cf. \cite{MRY}, \cite{RM} etc.), for instance, the vertex representations of the toroidal Lie algebras $T(\mathfrak g)$ were constructed by Moody-Rao-Yokonuma \cite{MRY} in simply laced types, generalizing the Frenkel-Kac construction
\cite{FK}.

As generalization of the quantum affine algebras, quantum toroidal algebras were introduced \cite{GKV} in the study of
geometric Langlands conjecture to algebraic surfaces. As quantum affinization, representations and structures of the quantum affine algebras
have been studied first for finite-dimensional representations (cf.\cite{CP1},\cite{CP3},\cite{H1} etc.) and then the vertex representations (cf. \cite{FJ}, \cite{J1}, \cite{J3}, \cite{J4}, \cite{JKK}, \cite{JKM}, \cite{JM} etc.). Varagnolo and Vasserot \cite{VV} proved a Schur-Weyl duality between the quantum toroidal algebras  $U_q({\mathfrak g}_{tor})$ and Cherednik's double affine Hecke algebras.
 The vertex representations of the quantum toroidal algebras in type $A$ were constructed \cite{S} as analogues of the Frenkel-Jing construction
\cite{FJ} of the quantum affine algebras, and then
the basic (vertex) representations of the quantum toroidal algebras in the simply laced
types were constructed \cite{FJW} via the McKay correspondence. In a series of papers \cite{M1, M2, M3}, Miki studied the structures
of the quantum toroidal algebra $U_q({\mathfrak g}_{tor})$ exclusively in type A.
Vertex operator realization of the quantum toroidal algebra $U_q({\mathfrak g}_{tor})$ of type $A$ was given in terms of
 the basic module for the affine Lie algebra $\hat{\mathfrak{gl}}_N$ \cite{GJ}. Much in depth work was
done for the quantum toroidal algebras of type $A$ (including the two parametric deformation) \cite{FJM1} \cite{FJM2},
 Yangians of toroidal types \cite{GTL} \cite{GM} and also in the viewpoint of quantum affinization \cite{H2}.
In spite of all these, it is still some away from full understanding of the quantum toroidal algebras in type A,
and even less is known for
the representation theory of other types.

Recently, quantum $N$-toroidal algebras
have been introduced by Gao, Jing, Xia, and Zhang \cite{GJXZ}
as natural generalization of the quantum toroidal algebras by replacing the $2$-torus with the $N$-torus.
The new quantum algebras can be realized as quotient algebras of the quantum GIM algebras \cite{Ko},
which generalize Slodowy's GIM algebras in connection with algebraic geometry \cite{Sl, BM}.
This offers a new perspective to understand quantum toroidal algebras. As nontrivial realization,
vertex representations of the quantum $N$-toroidal algebras
for simply-laced types have been constructed in \cite{GJXZ}, generalizing \cite{FJ}.

The goal of this paper is to construct a level-one module for the quantum $N$-toroidal algebra of type $C_n$, which provides new example of realization of the algebra for non-simply laced types and a model to understand enlargement of the torus.
Unlike the simply laced types \cite{FJ}, the vertex operators coming from the internal quantum Heisenberg subalgebra do not satisfy the
Drinfeld commutation relations in the non-simply laced types.  As in quantum affine algebras \cite{JKM}, we need to incorporate an
auxiliary bosonic field to the vertex
operator associated to each simple short root of the underlying
Lie algebra. In our construction, one sees that $N-1$ copies of the affine Heisenberg algebra stand at the common underlying lattice
to represent additional dimension in the torus and similarly the auxiliary Heisenberg subalgebra also bears similar property.
In the proof, we try to
offer more details on those special for non-simply laced situations.

The paper is organized as follows. In section 2, we recall the definition of the quantum $N$-toroidal algebra for type $C_n$ via generating functions.
We construct the Fock space and vertex operators in section 3. In the last section, we state our main result of the construction
and prove in detail that the module realizes the quantum $N$-toroidal algebra in type $C_n$ at level one.

\section{Quantum $N$-toroidal algebra $ U_q(\frak{g}_{N, tor})$}

In the paper, we always assume that $\frak{g}$ is the simple Lie algebra of type $C_{n}$ $(n\geq 3)$.
We review the notion of the quantum $N$-toroidal algebra $ U_q(\frak{g}_{N, tor})$ of type $C_n$
in this section.
For this, we start with some notations and data of the Lie algebra $\frak{g}$. Let $I=\{0,\cdots,n\}$, $I_0=\{1,\cdots,n\}$ and
$I_a=\{1, \ldots, n-1\}$.

Set $\mathbb{R}^{n}$ as the real Euclidean space with the inner product $(\hspace{0.2cm}|\hspace{0.2cm})$
and fix an orthogonal basis $\{e_1,\ldots,e_{n}\}$ such that $(e_i|e_j)=\frac{1}{2}\delta_{ij}$. The root system of $\frak{g}$ is spanned by the simple roots $\alpha_i$, where
$\alpha_i=e_i-e_{i+1}(i=1,\ldots,n-1), \alpha_n=2e_n$. So
$(\alpha_i|\alpha_i)=1$ for $1\leq i\leq n-1$ and $(\alpha_n|\alpha_n)=2$. Then the weight lattice and root lattice
of $\frak{g}$ are respectively
\begin{eqnarray*}
\dot{P}&=\{\sum_{i=1}^{n}k_ie_i|k_i \in \mathbb{Z}\},\\
\dot{Q}&=\{\sum_{i=1}^{n}k_i\alpha_i|k_i \in \mathbb{Z}\}.
\end{eqnarray*}

Denote by $\hat{\frak{g}}$  the affine Lie algebra associated to  $\frak{g}$, let $\delta$ be the null root of  $\hat{\frak{g}}$ satisfying $(\alpha_i|\delta)=0=(\delta|\delta),i \in I_{0}$, where the invariant bilinear form $(\ \ |\ \ )$ on $\hat{\frak{g}}$
extends that of $\frak{g}$. The set of simple roots
for the root system of
the affine Lie algebra $\hat{\frak{g}}$ (in type $C_{n}^{(1)}$) is given by 
$$\Pi=\{\alpha_0,\alpha_1,\cdots,\alpha_n\},$$
where $\alpha_0=\delta-2\alpha_1-2\alpha_2-\cdots-2\alpha_{n-1}-\alpha_n.$
The affine root and weight lattices are
\begin{eqnarray*}
Q=\sum_{i=0}^n\mathbb Z\alpha_i, \qquad\qquad P=\sum_{i=0}^n\mathbb Z\Lambda_i,
\end{eqnarray*}
where $\Lambda_i$ are defined by $\Lambda_i(\alpha_j^{\vee})=\delta_{ij}$ and $\Lambda_i(d)=0$ for the gradation $d$.
Here $\alpha_j^\vee=\frac{2\alpha_j}{(\alpha_j,\alpha_j)}$ for $j\in I$. Then the affine Cartan matrix $A=(a_{ij})$ is given by  $a_{ij}=(\alpha_j|\alpha_i^\vee)$ for $i,j\in I$.
So $d_ia_{ij}=(\alpha_i|\alpha_j)=d_ja_{ji}$ and 
$A$ is symmetrized by the diagonal matrix $D=diag(d_i|i\in I)$, where $(d_0, d_1,\cdots, d_n)=(1, \frac{1}{2},\cdots, \frac{1}{2}, 1)$
for type $C_n^{(1)}$.

As Lie algebra, the affine Lie algebra $\hat{\frak{g}}$ takes the form of the central extension of the loop algebra
${\mathfrak g}\otimes \mathbb C[t^{\pm1}]$, 
that is,  $\hat{\mathfrak{g}}=\mathfrak{g}\otimes\mathbb{C}[t,t^{-1}]\oplus\mathbb{C}\emph{\textbf{c}}$ with the
bracket $[x\otimes t^a,y\otimes t^b]=[x,y]\otimes t^{a+b}+a(x|y)\delta_{a+b,0}\emph{\textbf{c}}$, where $x,y\in \mathfrak g,a,b\in \mathbb{Z}$ and $\emph{\textbf{c}}$ is a nonzero central element.

For generic $q$, we set $q_i=q^{d_i},i\in I$.
Fix $N\geq 2$, let $J=\{1, \cdots, {N-1}\}$ and $e_s=(0, \cdots, 0, 1, 0, \cdots, 0)$ the $s$-th
standard unit vector of the $ (N-1)$-dimensional lattice $\mathbb Z^{N-1}$.  The general vector
of $\mathbb Z^{N-1}$ will be denoted by $\underline{k}=(k_1, k_2, \cdots, k_{N-1})$, and
$\underline{\bf{0}}$ is the ${(N-1)}$-dimensional zero vector.
For the formal variable $\underline{z}=(z_1, \cdots,\, z_{N-1})$, denote $\underline{z}^{\underline{k}}=\prod\limits_{s=1}^{N-1} z_s^{k_s}$.

We also need formal generating functions in the variables $\underline{z}$ and $z$, for instance, we set for each $s\in J$
\begin{gather*}
\delta(z)=\sum_{k\in\mathbb{Z}}z^{k},\qquad
  x_{i}^{\pm}(\underline{z})=\sum_{\underline{k} \in \mathbb{Z}^{N-1}}x_{i}^{\pm}(\underline{k}) \underline{z}^{-\underline{k}}, \qquad
  x_{i,s}^{\pm}({z_s})=\sum_{k_s\in \mathbb{Z}}x_{i}^{\pm}(k_se_s) z_s^{-k_s},\\
\phi_i^{(s)}(z) =\sum_{m \in \mathbb{Z}_+}\phi_i^{(s)}(m) z^{-m}, \qquad
\varphi_i^{(s)}(z)  = \sum_{n \in \mathbb{Z}_+}\varphi_i^{(s)}(-n) z^{n}.
\end{gather*}

Two versions of quantum numbers will be used in the paper. For a real number $a$, the $q$-analog of $a$ is defined by $(a;q)_{\infty}=\prod_{k=0}^{\infty}(1-aq^k)$. The quantum integer is defined as $[n]_q=\frac{q^n-q^{-n}}{q-q^{-1}}$
and the $q$-binomial numbers are $\Big[{m\atop  n}\Big]=\frac{[m]!}{[m-n]![n]!}$.

\begin{defi} \cite{GJXZ}
The quantum $N$-toroidal  algebra $ U_q(\mathfrak{g}_{N, tor})$
is the associative algebra with unit
$1$ generated by ${\gamma_s}^{\pm\frac{1}{2}}$, $s\in J$ and the coefficients of the generating functions 
$
x_i^{\pm}(\underline{z})$, $\phi_i^{(s)}(z_s)$, $\varphi_i^{(s)}(z_s)$, $i\in I
$
satisfying the relations below:
\begin{eqnarray}
&&\gamma_s^{\pm\frac1{2}}~~~ \hbox{are
central  such that } \gamma_s^{\frac{1}{2}}\gamma_s^{-\frac{1}{2}}=1, \\
&&\phi_{i}^{(s)}(0)\varphi_{j}^{(s')}(0)=\varphi_{j}^{(s')}(0)\phi_{i}^{(s)}(0)=1,\\
&& \phi_i^{(s)}(z_s)\phi_j^{(s')}(w_{s'})=\phi_j^{(s')}(w_{s'})\phi_i^{(s)}(z_s),
\quad\varphi_i^{(s)}(z_s)\varphi_j^{(s')}(w_{s'})=\varphi_j^{(s')}(w_{s'})\varphi_i^{(s)}(z),\\
&&g_{ij}\Bigl(z_s\ w_s^{-1}\gamma_s\Bigr)^{\delta_{s, s'}}\varphi_i^{(s)}(z_s)\phi_j^{(s')}(w_{s'})
=g_{ij}\Bigl(z_s\ w_s^{-1}\gamma_s^{-1}\Bigr)^{\delta_{s, s'}}
\phi_j^{(s')}(w_{s'})\varphi_i^{(s)}(z_s),\\
&&\varphi_i^{(s)}(z_s)x_j^{\pm}(\underline{w})\varphi_i^{(s)}(z_s)^{-1}
=g_{ij}\Bigl(\frac{z_s}{w_s}\gamma_s^{\mp
\frac{1}{2}}\Bigr)^{\pm1}x_j^{\pm}(\underline{w}),  \\
&&\phi_i^{(s)}(z_s)x_j^{\pm}(\underline{w})\phi_i^{(s)}(z_s)^{-1}=g_{ij}\Bigl(\frac{w_s}{z_s}\gamma_s^{\mp
\frac{1}{2}}\Bigr)^{\mp1}x_j^{\pm}(\underline{w}), \\
&&\lim_{z\to w}[x_{i,s}^{\pm}(z), x_{i,s'}^{\pm}(w)]=0, \qquad \hbox{for} \quad s\neq s',\label{e:comm06}
\end{eqnarray}
\begin{eqnarray}
&&(z-q^{\pm (\alpha_i|\alpha_j)}w)\,x_{i, s}^{\pm}(z)x_{j, s}^{\pm}(w)=(q^{\pm (\alpha_i|\alpha_j)}z-w)\,x_{j, s}^{\pm}(w)\,x_{i, s}^{\pm}(z)\label{e:comm07},\\
&&[\,x_{i, s}^{+}(z),x_{j, s}^{-}(w)\,]=\frac{\delta_{ij}}{q_i-q^{-1}_i}\Big(\phi_i^{(s)}(w\gamma_s^{\frac{1}2})\delta(\frac{w\gamma_s}{z})
-\varphi_i^{(s)}(w\gamma_s^{-\frac{1}{2}})\delta(\frac{w\gamma^{-1}_s}{z})\Big),\label{e:comm08}
\end{eqnarray}
where $\phi_i^{(s)}(r)$ and $\varphi_i^{(s)}(-r)\, (r\geq 0)$ such that $\phi_i^{(s)}(0)=K_i$ and  $\varphi_i^{(s)}(0)=K_i^{-1}$
are defined by:
\begin{gather*}\sum\limits_{r=0}^{\infty}\phi_i^{(s)}(r) z^{-r}=K_i \exp \Big(
(q_i{-}q_i^{-1})\sum\limits_{\ell=1}^{\infty}
a_i^{(s)}(\ell)z^{-\ell}\Big), \\
\sum\limits_{r=0}^{\infty}\varphi_i^{(s)}(-r) z^{r}=K_i^{-1}\exp
\Big({-}(q_i{-}q_i^{-1})
\sum\limits_{\ell=1}^{\infty}a_i^{(s)}(-\ell)z^{\ell}\Big),
\end{gather*}
\begin{eqnarray}\label{e:comm09}
&&\underset{z_1,\ldots,z_m}{Sym}\sum_{k=0}^{m=1-a_{ij}}(-1)^k\Big[{m\atop  k}\Big]_{i}x_{i, s}^{\pm}({z_1})\cdots x_{i, s}^{\pm}(z_k) x_{j,s}^{\pm}(w)x_{i, s}^{\pm}({z_{k+1}})    
x_{i, s}^{\pm}({z_m})=0, \\
&&\lim_{z_k\to w} \sum_{k=0}^{3}(-1)^k
	\Big[{3\atop  k}\Big]_{i}x_{i,s}^{\pm}(z_1))\cdots x_{i,s}^{\pm}(z_k)x_{i,s'}^{\mp}(w) x_{i,s}^{\pm}(z_{k+1})\cdots x_{i,s}^{\pm}(z_3)=0,
	~~~s\neq s'\in J,\label{e:comm10}
\end{eqnarray}
where $g_{ij}^{\pm}(z):=\sum_{n\in
\mathbb{Z}_+}c^{\pm}_{ijn}z^{n}$ is the Taylor series expansion of $g_{ij}^{\pm}(z)=(\frac{zq_i^{a_{ij}}-1}{z-q_i^{a_{ij}}})^{\pm 1}$ at $z=0$
in $\mathbb{C}$.
\end{defi}
Alternatively, the quantum algebra has the following equivalent definition \cite{GJXZ}.
\begin{defi}
\,The quantum $N$-toroidal algebra $ U_q(\frak{g}_{N, tor})$ is a complex unital associative algebra with
generators $x_{i}^{\pm}(\underline{k}),\, a_i^{(s)}(r),\, K_i^{\pm}$ and $\gamma_s^{\pm\frac{1}{2}}$  $(i\in I,\, s\in J,\,\underline{k}\in\mathbb{Z}^{N-1},\, r\in \mathbb{Z}/\{0\})$ and the following relations:
\begin{eqnarray}
&&\gamma_s^{\pm\frac{1}2} \textrm{are central,}\  K_i^{\pm1}~~~~ \textrm{and}~~~~  a_j^{(s)}(r) ~~~~\textrm{ commute with each other}\label{comm02},\\
&&[\,a_i^{(s)}(r),a_j^{(s')}(l)\,]
=\delta_{s,s'}\delta_{r+l,0}\frac{[\,r\,a_{ij}\,]_i}{r}
\frac{\gamma_s^{r}-\gamma_s^{-r}}{q_j-q_j^{-1}}\label{comm03},\\
&&K_ix_{j}^{\pm}(\underline{k})K_i^{-1}=q^{\pm a_{ij}}_ix_{j}^{\pm}(\underline{k})\label{comm04},\\
&&[\,a_i^{(s)}(r),x_{j}^{\pm}(\underline{k})\,]=\pm \frac{[\,r\,a_{ij}\,]_i}{r}
\gamma_s^{\mp\frac{|r|}{2}}x_{j}^{\pm}(re_s{+}\underline{k})\label{comm05},\\
&&[x_i^{\pm}(ke_s),x_{i}^{\pm}(le_{s'})\,]=0, \quad \hbox{for}\quad s\neq s' \quad \hbox{and}\quad kl\neq 0, \label{comm06}\\
&&[\,x_{i}^{\pm}((k+1)e_s),\,x_{j}^{\pm}(l_se_s)\,]_{q^{\pm (\alpha_i, \alpha_j)}}+
[\,x_{j}^{\pm}((l_s+1)e_s),\,x_{i}^{\pm}(ke_s)\,]_{q^{\pm (\alpha_i, \alpha_j)}}=0,\label{comm07}\\
&&[\,x_{i}^{+}(ke_{s}),\,x_j^{-}(le_{s})\,]=\delta_{ij}\big(\frac{\gamma_s^{\frac{k-l}{2}}\phi_i^{(s)}((k+l))
-\gamma_s^{\frac{l-k}{2}}\varphi_i^{(s)}((k+l))}{q_{i}-q_{i}^{-1}}\big),\label{comm08}
\end{eqnarray}
where $\phi_i^{(s)}(r)$ and $\varphi_i^{(s)}(-r)\, (r\geq 0)$ 
are defined by
\begin{gather*}\sum\limits_{r=0}^{\infty}\phi_i^{(s)}(r) z^{-r}_s=K_i \exp \Big(
(q_i{-}q_i^{-1})\sum\limits_{\ell=1}^{\infty}
a_i^{(s)}(\ell)z^{-\ell}\Big), \\
\sum\limits_{r=0}^{\infty}\varphi_i^{(s)}(-r) z^{r}_s=K_i^{-1}\exp
\Big({-}(q_i{-}q_i^{-1})
\sum\limits_{\ell=1}^{\infty}a_i^{(s)}(-\ell)z^{\ell}\Big),
\end{gather*}
\begin{eqnarray}\label{comm09}
 && \underset{k_1, \ldots, k_m}{Sym}\sum_{t=0}^{m=1-a_{ij}}(-1)^t\Big[{m\atop  t}\Big]_{i}x_i^{\pm}(k_te_s)\cdots
x_i^{\pm}(k_te_s)x_{j}^{\pm}(\ell e_s)x_i^{\pm}(k_{t+1}e_s)\cdots x_i^{\pm}(k_me_s)=0,
\\\label{comm10}
&&\sum_{k=0}^{3}(-1)^k
	\Big[{3\atop  k}\Big]_{i}x_i^{\pm}({e_sm_1})\cdots x_i^{\pm}({e_sm_{k}}) x_{i}^{\mp}({e_{s'}\ell})x_i^{\pm}({e_sm_{k+1}})\cdots x_i^{\pm}({e_sm_{3}})=0,\\ \nonumber
	&&\hspace{3.1cm} \qquad	~~~~\hbox{for}~~~~  i\in I_0~~~~ \hbox{and}~~~~m_1m_2m_3\ell\neq 0, \, s\neq s'\in J.
\end{eqnarray}
Here $[m]_i=[m]_{q_i}$
and the $q$-bracket is defined by $[a, b]_{u}\doteq ab-uba$.
\end{defi}

\begin{remark}\, Relation \eqref{e:comm09} or \eqref{comm09} is void for $i=j$, as $m=-1$ in this case.
The defining relations $\eqref{comm07}$ to $\eqref{comm10}$ are written in terms of the components $x_i^{\pm}(ke_s)$.
Using relation $\eqref{comm05}$, we can obtain the relations for the generators $x_i^{\pm}(\underline{k})$.
\end{remark}

\begin{remark}
By relations $\eqref{comm03}$, $\eqref{comm05}$ and $\eqref{comm08}$, we have more identities as follows. 
It follows  from $\eqref{comm08}$ that
$$a_i^{(s)}(1)=K_i^{-1}\gamma_s^{1/2}\,[\,x_i^+(\underline{0}),\,x_i^-(e_s)\,],\qquad
a_i^{(s')}(-1)=K_j{\gamma_{s'}}^{-1/2}\,[\,x_i^+(-e_{s'}),\,x_i^-(\underline{0})\,].$$

Using the first relation to expand the bracket $[a_i^{(s)}(1), a_j^{(s')}(-1)]$ 
we obtain for $s\neq s$
\begin{equation*}
\begin{split}
0&=[\,a_i^{(s)}(1), a_j^{(s')}(-1)\,]=K_i^{-1}\gamma_s^{\frac{1}{2}}[\,[x_i^+(\underline{0}),
x_i^-(e_s)], a_j^{(s')}(-1)\,]\\
&=K_i^{-1}\gamma_s^{\frac{1}{2}}\Bigl([\,[x_i^+(\underline{0}), a_j^{(s')}(-1)],
x_i^-(e_s)\,]+[\,x_i^+(\underline{0}), [x_i^-(e_s), a_j^{(s')}(-1)]\,]\Bigr)\\
&=[-a_{ij}]_iK_i^{-1}\gamma_s^{\frac{1}{2}}\Bigl(\gamma_{s'}^{-\frac{1}{2}}[x_i^+(-e_{s'}),
x_i^-(e_s)]-\gamma_{s'}^{\frac{1}{2}}[x_i^+(\underline{0}),x_i^-(e_s-e_{s'})]\Bigr).
\end{split}
\end{equation*}

Therefore, 
$$[x_i^+(\underline{0}),x_i^-(e_s-e_{s'})]=\gamma_s^{-1}[x_i^+(-e_{s'}),
x_i^-(e_s)].$$

Similarly expanding the bracket by using the second relation, we also get that 
\begin{equation*}
\begin{split}
0&=[\,a_i^{(s)}(1), a_j^{(s')}(-1)\,]=K_j\gamma_{s'}^{-\frac{1}{2}}[\,a_i^{(s)}(1), [x_j^+(-e_{s'}),
x_j^-(\underline{0})] \,]\\
&=K_j\gamma_{s'}^{-\frac{1}{2}}\Bigl([\,[a_i^{(s)}(1), x_j^+(-e_{s'})], x_j^-(\underline{0})\,]
+[\,x_j^+(-e_{s'}), [a_i^{(s)}(1), x_j^-(\underline{0})]\,]\Bigr)\\
&=[a_{ij}]_iK_j\gamma_{s'}^{-\frac{1}{2}}\Bigl(\gamma_{s}^{-\frac{1}{2}}[x_j^+(e_s-e_{s'}),
x_j^-(\underline{0})]-\gamma_{s}^{\frac{1}{2}}[x_j^+(-e_{s'}),x_j^-(e_s)]\Bigr).
\end{split}
\end{equation*}

So, 
$$[x_j^+(e_s-e_{s'}),
x_j^-(\underline{0})]=\gamma_s^{-1}[x_j^+(-e_{s'}),x_j^-(e_s)].$$

\end{remark}

\section{Fock space and vertex operators}

We construct the Fock space and vertex operators for the quantum $N$-toroidal algebra $ U_q(\frak{g}_{N, tor})$ for type $C_n$
in this section.

\subsection{Fock space}
Let $\tilde{Q}$ be an identical copy of the root lattice of type $A_{n-1}$, or the lattice of short roots in $\mathfrak g$,
so $\tilde{Q}=\oplus_{i=1}^{n-1}\mathbb Z\tilde{\alpha}_i$. We equip $Q[A_{n-1}]$ with the inner product
$(\tilde{\alpha}_i|\tilde{\alpha}_j)=(\alpha_i|\alpha_j)=\delta_{ij}-\frac{1}{2}\delta_{|i-j|,1}.$

We also adopt the nondegerate bilinear form $(\ |\ )$ on the lattice $\mathbb ZJ$ defined by
\begin{equation}
(i|j)=\begin{cases} -1 & i\neq j\\ 0 & i=j \end{cases},
\end{equation}
where the elements of $J$ are written as $\{1, 2, \cdots, {N-1}\}$.

Let $\mathcal{H}$ be the quantum Heisenberg Lie algebra generated by  $a^{(s)}_i(n)$ $(n\in \mathbb{Z}^{\times} ,i\in I, s\in J),a_i(0)$  and $b^{(s)}_i(n), b_i(0)$ $(n\in \mathbb{Z}^{\times}, i\in I_a, s\in J)$  satisfying:
   \begin{eqnarray*}
   &&[a^{(s)}_i(m),a^{(s^{'})}_j(l)]=\delta_{s,s^{'}}\delta_{m,-l}\frac{[a_{ij}m]_i}{m}\frac{q^m-q^{-m}}{q_j-q^{-1}_j},\\ &&[b^{(s)}_i(m),b^{(s^{'})}_j(l)]=\delta_{s,s^{'}}\delta_{m,-l}\frac{[(\tilde{\alpha}_i|\tilde{\alpha}_j)m]}{m}[m],\\
   &&[a^{(s)}_i(m),b^{(s')}_j(l)]=0,\\
   &&[a_i(0),a^{(s^{'})}_j(m)]=[b_i(0),b^{(s^{'})}_j(m)]=0,
   \end{eqnarray*}
where $m, l\in\mathbb Z^{\times}, s, s'\in J$ and $i, j\in I$ or $I_a$ depending on $a^{(s)}_i(m)$ or $b^{(s)}_i(m)$.

Let $S(\mathcal{H}^-)$ be the symmetric algebra generated by $\{a_i^{(s)}(-m),b_j^{(s)}(-m)|m\in \mathbb{N},i\in I, j\in I_a, s\in J\}$. We define the Fock space
$$\mathcal{V}=S(\mathcal{H}^-)\otimes \mathbb{C}[Q]\otimes \mathbb{C}[\tilde{Q}]\otimes \mathbb C[\mathbb ZJ],$$
where $\mathbb{C}[\tilde{Q}]=\mathbb C[\underset{i\in I_a}{\oplus}\mathbb{Z}\tilde{a}_i]$ is the group algebra of $\tilde{Q}$, and similarly for
$\mathbb C[Q]$ and $\mathbb C[\mathbb ZJ]$.

The space $\mathcal V$ becomes a natural $\mathcal H$-module by
letting $a_i^{(s)}(-m), b_i^{(s)}(-m)$ $(m>0)$ act as multiplication operators, and
$a_i^{(s)}(m), b_i^{(s)}(m)$ $(m>0)$ as differentiation operators subject to the Heisenberg relations.

Define the operators $a_i(0),b_j(0), s(0)$ ($i\in I, j\in I_a, s\in J$)
and  $e^{\alpha}, e^{\tilde{\alpha}}, e^s$ ($\alpha\in\mathbb{C}[Q], \tilde{\alpha}\in \mathbb{C}[\tilde{Q}], s\in \mathbb ZJ$) act on $\mathcal{V}$ by
 \begin{gather*}
   {a_i(0)e^{\lambda}e^{\tilde{\lambda}}e^{t}=(\alpha_i|\lambda)e^{\lambda}e^{\tilde{\lambda}}}e^t,
   \quad
   b_j(0)e^{\lambda}e^{\tilde{\lambda}}e^t
   =(\tilde{\alpha}_j|\tilde{\lambda})e^{\lambda}e^{\tilde{\lambda}}e^t, \quad s(0)e^{\lambda}e^{\tilde{\lambda}}e^t=(s|t)e^{\lambda}e^{\tilde{\lambda}}e^t,\\
   e^{\alpha}e^{\lambda}e^{\tilde{\lambda}}=e^{\alpha+\lambda}e^{\tilde{\lambda}}e^t,
   \quad e^{\tilde{\alpha}}e^{\lambda}e^{\tilde{\lambda}}=e^{\lambda}e^{\tilde{\alpha}+\tilde{\lambda}}e^t, \quad
   e^{s}e^{\lambda}e^{\tilde{\lambda}}=e^{\alpha}e^{\tilde{\lambda}}e^{s+t},
   \end{gather*}
   {where $\lambda\in \mathbb{C}[Q]$, $\tilde{\lambda}\in \mathbb{C}[\tilde{Q}], t\in \mathbb{C}[\mathbb ZJ]$.

Let $z$ be a complex variable, the above implies that
$$z^{a_i(0)}e^{\lambda}e^{\tilde{\lambda}}e^t=z^{(\alpha_i|\lambda)}e^{\lambda}e^{\tilde{\lambda}}e^t,
\quad z^{b_j(0)}e^{\lambda}e^{\tilde{\lambda}}e^t
=z^{(\tilde{\alpha}_j|\tilde{\lambda})}e^{\lambda} e^{\tilde{\lambda}}e^t, \quad z^{s(0)}e^{\lambda}e^{\tilde{\lambda}}e^t=z^{(s|t)}e^{\lambda}e^{\tilde{\lambda}}e^t. $$

We introduce the sign operators $\varepsilon_i$ on $\mathcal V$ such that 
\begin{eqnarray*}
&&\varepsilon_i\varepsilon_j=(-1)^{2(\alpha_i|\alpha_j)}\varepsilon_j\varepsilon_i, \ \ \ \ \ \ \text{for all}\ \ 1\leq i,j\leq n-1,\ \ \ \ \ \ \ \ \ \ \ \ \ \ \ \ \ \ \\
&&\varepsilon_i\varepsilon_j=(-1)^{(\alpha_i|\alpha_j)}\varepsilon_i\varepsilon_j, \ \ \ \ \ \ \ \ \ \ i \ \text{or}\  j=0,n,\\
&& (\varepsilon_i)^2=1, \qquad\forall\, i.
\end{eqnarray*}

We remark that these operators can be constructed from a cocycle on the abelian group $\mathbb Z_2^{n+1}$.

\subsection{Normal order} The normal order :\hspace{0.4cm}: is defined as usual:
\begin{eqnarray*}
:a_i^{(s)}(m)a_j^{(s)}(n):&=&\left\{\begin{array}{rcl}a_i^{(s)}(m)a^{(s)}_j(n),&&{m\leq n},\\a^{(s)}_i(n)a^{(s)}_j(m),&& {m>n},\end{array} \right.\\
{:e^{\alpha}a_i(0):} &=& :a_i(0)e^{\alpha}: = e^{\alpha}a_i(0),\\
:e^{\tilde{\alpha}}b_i(0): &=& :b_i(0)e^{\tilde{\alpha}}: = e^{\tilde{\alpha}}b_i(0),\\
:e^ss'(0):&=&=:s'(0)e^s:=e^ss'(0),
\end{eqnarray*}
and similarly for $b^{(s)}_i(m)$. For commuting operators $a_i^{(s)}(m), a_i^{(s')}(m)$ the normal order product is the usual product.

It is easy to see that $a_i(0),b_j(0), s(0),e^{\alpha_j}, e^{\tilde{\alpha}_k}, e^{s'}$ commute with each other except that
$$[a_i(0),e^{\alpha_j}]=(\alpha_i|\alpha_j)e^{\alpha_j},\quad [b_j(0),e^{\tilde{\alpha}_k}]=(\tilde{\alpha}_j|\tilde{\alpha}_k)e^{\tilde{\alpha}_k},
\quad [s(0), e^{s'}]=(\delta_{ss'}-1)e^{s'}.$$

\subsection{Vertex operators}
Following \cite{JKM}, we introduce the main vertex operators of $ U_q(\frak{g}_{N, tor})$ as follows.
First we define the normal order product $:\ \ :$ of the vertex operators by moving the annihilation operators like $a_i^{(s)}(k), b_j^{(s)}(k)$ $(k>0)$ etc. to the right and keep the sign operators $\varepsilon_i$ in the order inside (for $:XY:$ and like).

For
$i\in I, j\in I_a, s\in J$ and $\epsilon=\pm$, we define
   \begin{eqnarray*}
  &&Y^{\pm}_{i,s}(z)=\exp\Big(\pm\sum\limits_{k=1}^{\infty}\frac{a^{(s)}_i(-k)}{[k/d_i]_i}q^{\mp\frac{k}{2}}z^{k}\Big)
   \exp\Big(\mp\sum\limits_{k=1}^\infty\frac{a^{(s)}_i(k)}{[k/d_i]_i}q^{\mp\frac{k}{2}}z^{-k}\Big){e^{\pm{\alpha_i}}z^{\pm a_i(0)+d_i}\varepsilon_i}\\
  &&{\hskip 1.1cm}=:\exp\Big(\mp\sum\limits_{k\neq 0}^\infty\frac{a^{(s)}_i(k)}{[k/d_i]_i}q^{\mp\frac{|k|}{2}}z^{-k}\Big):e^{\pm{\alpha_i}}
  z^{\pm a_i(0)+d_i}\varepsilon_i,\\
 &&Y^{\pm}_{i}(\underline{z})=
   :\exp\Big(\mp\sum\limits_{s=1}^{N-1}\sum\limits_{k_s\neq 0}^\infty\frac{a^{(s)}_i(k_s)}{[k_s/d_i]_i}q^{\mp\frac{|k_s|}{2}}z_s^{-k_s}\Big):
   {e^{\pm{\alpha_i}}\prod_{s\in J}z_s^{\pm a_i(0)+d_i}\varepsilon_i},
 \end{eqnarray*}
 \begin{eqnarray*}
   &&{\hskip -2.8cm}X_{j\epsilon,s}^{\pm}(z)=:\exp\Bigg(\sum\limits_{k\neq 0}^{\infty}\Big(\mp\frac{q^{\mp\frac{|k|}2}z^{-k}}{[k/d_j]_j}a^{(s)}_j(k)-\epsilon \frac{(q^{\frac{\epsilon}{2}}z)^{-k}}{[k]}b^{(s)}_j(k)\Big)\Bigg):\\
   &&\cdot e^{\pm{\alpha_j}+\epsilon\tilde{\alpha}_j\pm s}z^{\pm a_j(0)+\epsilon b_j(0)+1}q^{\frac{ b_j(0)}{2}\pm s(0)}(-1)^{(1\mp\epsilon)b_j(0)}\varepsilon_j,
 \end{eqnarray*}
 \begin{eqnarray*}
   &&X_{j\epsilon}^{\pm}(\underline{z})=
   :\exp\Bigg(\sum\limits_{s=1}^{N-1}\sum\limits_{k_s\neq 0}^{\infty}\Big(\mp\frac{q^{\mp\frac{|k_s|}2}z^{-k_s}}{[k_s/d_j]_j}a^{(s)}_j(k_s)-\epsilon \frac{(q^{\frac{\epsilon}{2}}z_s)^{-k_s}}{[k_s]}b^{(s)}_j(k_s)\Big)\Bigg):\\
   &&\hskip1.3cm\cdot e^{\pm{\alpha_j}+\epsilon\tilde{\alpha}_j\pm s}\prod\limits_{s=1}^{N-1}z_s^{\pm a_j(0)+\epsilon b_j(0)+1}q^{\frac{ b_j(0)}{2}\pm s(0)}(-1)^{(1\mp\epsilon)b_j(0)}\varepsilon_j^{(s)}, 
  \end{eqnarray*}
 \begin{eqnarray*}
  &&\Phi_i^{(s)}(z)=q^{a_i(0)}\exp \Big((q_i{-}q_i^{-1})\sum\limits_{\ell=1}^{\infty}a_i^{(s)}(\ell)z^{-\ell}\Big), \\
   &&\Psi_i^{(s)}(z)=q^{-a_i(0)}\exp\Big({-}(q_i{-}q_i^{-1})\sum\limits_{\ell=1}^{\infty}a_i^{(s)}(-\ell)z^{\ell}\Big).
   \end{eqnarray*}
Note that we have used $\varepsilon_i^{(s)}$ to mean the action of $\varepsilon_i$ on the subspace indexed by $s$.
Therefore
\begin{eqnarray*}
&&X^{\epsilon'}_{i\epsilon, s}(z)=
:\exp\Bigg(-\sum\limits_{k\neq 0}^{\infty}\Big(\epsilon'\frac{q^{{-\epsilon'}\frac{|k|}2}z^{-k}}{[k/d_j]_j}a^{(s)}_i(k)+\epsilon \frac{(q^{\frac{\epsilon}{2}}z)^{-k}}{[k]}b^{(s)}_i(k)\Big)\Bigg):\\
&&\qquad\qquad\qquad \cdot e^{\epsilon'{\alpha_i}+\epsilon\tilde{\alpha}_i+\epsilon's}z^{\epsilon' a_i(0)+\epsilon b_i(0)+1}q^{\frac{ b_j(0)}{2}+\epsilon's(0)}(-1)^{(1-\epsilon'\epsilon)b_i(0)}\varepsilon_i.
\end{eqnarray*}

\begin{prop}
The operator product expansions (OPE) are given  for $i,j\in I$, $k,l\in\{1,\cdots,n-1\}$, $\epsilon,\epsilon'\in \{\pm\}$ and $s\in J$ as follows.
 \begin{eqnarray}\notag
&& Y^{\epsilon}_{i,s}({z})Y^{\epsilon'}_{j,s}({w})=:Y^{\epsilon}_{i,s}({z})Y^{\epsilon'}_{j,s}({w}):\\ \label{3.1}
&&\left\{\begin{array}{lll}
 \ \ \ \ \ \ \ \ \ \ \ \ 1,&&{(\alpha_i|\alpha_j)=0},
\\(z-wq^{-(\epsilon+\epsilon')/2})^{\pm\epsilon\epsilon'},&&(\alpha_i|\alpha_j)=\pm 1,
\\(z-wq^{1-(\epsilon+\epsilon')/2})^{\epsilon\epsilon'}(z-wq^{-1-(\epsilon+\epsilon')/2})^{\epsilon\epsilon'},&&(\alpha_i|\alpha_j)=2,\\
f(z-wq^{-\frac{\epsilon+\epsilon'}2})^{\epsilon\epsilon'},
 &&(\alpha_i|\alpha_j)=-\frac{1}{2}.
\end{array} \right.
 \end{eqnarray}
\begin{eqnarray}\notag
&&X_{k\epsilon_1,s}^{\epsilon}(z)X_{l\epsilon_2,s}^{\epsilon'}(w)=:X_{k\epsilon_1,s}^{\epsilon}(z)X_{l\epsilon_2,s}^{\epsilon'}(w):\\ \label{3.2}
&&\left\{\begin{array}{lll}
\ \ \ \ \ \ \ \ \ \ \ \ \ \ \ \ \ \ \ \ \ \ \ \ \ \ \ \ \ \ \ \ \ \ \ \ 1,&&{(\alpha_k|\alpha_l)=0},
\\(z-q^{-(\epsilon+\epsilon')/2}w)^{\epsilon\epsilon'}(q^{\epsilon_1/2}z-q^{\epsilon_2/2}w)^{\epsilon_1\epsilon_2},&&(\alpha_k|\alpha_l)=1,
\\ f(z-q^{-\frac{\epsilon+\epsilon'}2}w)^{\epsilon\epsilon'}
f(q^{\frac{\epsilon_1}2}z-q^{\frac{\epsilon_2}2}w)^{\epsilon_1\epsilon_2}
(-1)^{(1-\epsilon\epsilon_1)/2}, && (\alpha_k|\alpha_l)=-\frac12,
\end{array} \right.    
 \end{eqnarray}
 {\begin{eqnarray}
&&X_{k\epsilon_1,s}^{\epsilon}(z)X_{k\epsilon_2,s'}^{\epsilon'}(w)
=:X_{k\epsilon_1,s}^{\epsilon}(z)X_{k\epsilon_2,s'}^{\epsilon'}(w): z^{\epsilon\epsilon'+\epsilon_1\epsilon_2},  
\qquad\hbox{for}\quad s\neq s',\label{3.3}
\end{eqnarray}}
where $f(a-x)=\frac{(xq^{3/2}/a; q^2)_{\infty}}{(xq^{1/2}/a; q^2)_{\infty}}a^{-1/2}$
a $q$-analog of the infinite series $(a-x)^{-1/2}$ at $x=0$ \cite{J3}.
\begin{eqnarray}\label{3.4}
&&Y_{i,s}^{\epsilon}(z)X_{k\epsilon,s}^{\epsilon'}(w)=:Y_{i,s}^{\epsilon}(z)X_{k\epsilon,s}^{\epsilon'}(w):
(z-q^{-(\epsilon+\epsilon')/2}w)^{-\epsilon\epsilon'},\qquad (\alpha_i|\alpha_k)=-1,\\
&&X_{k\epsilon,s}^{\epsilon}(z)Y_{i,s}^{\epsilon'}(w)=:X_{k\epsilon,s}^{\epsilon}(z)Y_{i,s}^{\epsilon'}(w):
(z-q^{-(\epsilon+\epsilon')/2}w)^{-\epsilon\epsilon'},\qquad (\alpha_i|\alpha_k)=-1.
\end{eqnarray}
\end{prop}

\begin{remark} The normal order product $:X_{k\epsilon_1,s}^{\epsilon}(z)X_{l\epsilon_2,s}^{\epsilon'}(w):$
contains the term $(-1)^{\frac{1-\epsilon\epsilon_1}2b_k(0)}$ $(-1)^{\frac{1-\epsilon'\epsilon_2}2b_l(0)}$.
Also $(-1)^{2b_k(0)}=1$.
\end{remark}

\begin{proof}
It is enough to check \eqref{3.1}, as \eqref{3.2} can be obtained similarly.
Five cases need to be considered: $(\alpha_i|\alpha_j)=0, 1, -1, 2 $ and $-\frac{1}{2}$.
Since all verifications are similar, we just give the case of $(\alpha_i|\alpha_j)=1$.
It is clear that
\begin{eqnarray*}
&&Y^{\epsilon}_{i,s}({z})Y^{\epsilon'}_{i,s}({w})=:Y^{\epsilon}_{i,s}({z})Y^{\epsilon'}_{i,s}({w}):
 exp\big(-\sum\limits_{n\geq1}\frac{[\epsilon\epsilon'n]}{n[n]}(\frac{w}{zq^{(\epsilon+\epsilon')/2}})^n\big)z^{\epsilon\epsilon'}\\
&&\hskip2.3cm     =:Y^{\epsilon}_{i,s}({z})Y^{\epsilon'}_{i,s}({w}):
(z-wq^{-(\epsilon+\epsilon')/2})^{\epsilon\epsilon'}.
\end{eqnarray*}
\end{proof}

\section{Vertex representation of $ U_q(\frak{g}_{N, tor})$}

The purpose of this section is to give the main result and the proof. 

\theoremstyle{theorem}
\begin{theorem}\label{theorem1}

The map $\pi$ from $ U_q(\frak{g}_{N, tor})$ to $End(\mathcal{V})$ defined by
$$\pi(K_i)={q^{a_i(0)}}, \pi(a^{(s)}_i(m))=a^{(s)}_i(m), \pi(\gamma_s)=q,$$
$$\pi(\phi_i^{(s)}(z))=\Phi_i^{(s)}(z), \pi(\varphi_i^{(s)}(z))=\Psi_i^{(s)}(z), $$
$$\pi(x^{\pm}_{j,s}({z}))=X^{\pm}_{j,s}({z})=\sum_{\epsilon\in\{\pm\}}X_{j\epsilon,s}^{\pm}(z),\pi(x^{\pm}_{j}({\underline{z}}))=X^{\pm}_{j}({\underline{z}})=\sum_{\epsilon\in\{\pm\}}X_{j\epsilon}^{\pm}(\underline{z}),$$
$$\pi(x^{\pm}_{k,s}({z}))=X^{\pm}_{k,s}({z})=Y^{\pm}_{k,s}(z),\pi(x^{\pm}_{k}({\underline{z}}))=X^{\pm}_{k}({\underline{z}})=Y^{\pm}_{k}(\underline{z}),$$
$i\in I,j\in I_a=\{1,\cdots,n-1\},k\in\{0,n\}$
 gives a representation of $ U_q(\frak{g}_{N, tor})$.
\end{theorem}

\begin{remark}
Since our construction already obeys the relations $(2.1)$-$(2.4)$, we only need to check the relations $(2.5)$-$(2.9)$.
\end{remark}

\subsection{Proof of relations $(2.5)$ and $(2.6)$.}
As relation $(2.6)$ is similar to relation $(2.5)$, it is enough to check relation $(2.5)$.
In fact, by moving $\Psi_i^{(s)}(z_s)$ across $Y_j^+(\underline{w})$ we have that 

\begin{eqnarray*}
&&\Psi_i^{(s)}(z_s)Y_{j}^{+}(\underline{w})\\
&=&{q^{-(\alpha_i|\alpha_j)}\exp(\sum\limits_{l>0}\frac{q_i^{-la_{ij}}-q_i^{la_{ij}}}{l}(\frac{q^{-\frac{1}{2}}z_s}{w_s})^l) Y_{j}^{+}(\underline{w})\Psi_i^{(s)}(z_s)}\\
&=&g_{ij}\Bigl(\frac{z_s}{w_s}q^{-
\frac{1}{2}}\Bigr)Y_{j}^{+}(\underline{w})\Psi_i^{(s)}(z_s).
\end{eqnarray*}
{\subsection{Proof of relation $(2.7)$.}
We just show relation $(2.7)$ for the short roots, i.e. $i\in I_a$. For $s\neq s'$, we have that
\begin{eqnarray*}
X^{+}_{i, s}(z)X^{+}_{i, s'}(w)=\sum_{\epsilon_1, \epsilon_2}:X^{+}_{i\epsilon_1, s}(z)X^{+}_{i\epsilon_2, s'}(w): z^{1+\epsilon_1\epsilon_2}. \end{eqnarray*}

 Therefore $\lim_{z\to w}[X_{i,s}^{+}(z), X_{i,s'}^{+}(w)]=0$. 

\subsection{Proof of relation $(2.8)$.}
There are five cases to consider: $(\alpha_i|\alpha_j)=0, -\frac{1}{2}, \pm 1, $ and $2$. The first and the last are trivial.

Now let $(\alpha_i|\alpha_j)=-1$, by \eqref{3.1}-\eqref{3.2} it follows that
$$
X_{n-1,s}^{\pm}({z})X_{n,s}^{\pm}({w})=:X_{n-1,s}^{\pm}({z})X_{n,s}^{\pm}({w}):(z-q^{\mp 1}w)^{-1}.
$$

Since $:X_{n-1,s}^{-}({z})X_{n,s}^{-}({w}):=-:X_{n,s}^{-}({w})X_{n-1,s}^{-}({z}):$, we get that
\begin{eqnarray*}
 (z-q^{\mp 1}w)\,X_{n-1,s}^{\pm}({z})X_{n,s}^{\pm}({w})
&=&(q^{\mp 1}z-w)\,X_{n,s}^{\pm}({w})\,X_{n-1,s}^{\pm}({z}).
\end{eqnarray*}

Take the case $(\alpha_i|\alpha_j)=-\frac{1}{2}$, $i, j\in I_a$. Note that $f(a-x)_{q\to q^{-1}}=f(a-x)$ and
\begin{eqnarray}\label{fcn1}
f(a-x)f(a-q^{-1}x)&=&\frac{(xq^{\frac{3}{2}}/a; q^2)_{\infty}}{(xq^{\frac 12}/a; q^2)_{\infty}}
\frac{(xq^{\frac12}/a; q^2)_{\infty}}{(xq^{-\frac12}/a; q^2)_{\infty}}a^{-1}=(a-q^{-\frac12}x)^{-1}, \\ \label{fcn2}
\frac{f(a-q^{-1}x)}{f(a-qx)}&=&\frac{(xq^{\frac12}/a; q^2)_{\infty}}{(xq^{-\frac12}/a; q^2)_{\infty}}
\frac{(xq^{\frac32}/a; q^2)_{\infty}}{(xq^{\frac52}/a; q^2)_{\infty}}=\frac{a-q^{\frac12}x}{a-q^{-\frac12}x}.
\end{eqnarray}

It follows from \eqref{3.1}-\eqref{3.2} that
\begin{eqnarray*}
X_{i,s}^{+}({z})X_{j,s}^{+}({w})=\sum_{\epsilon_1, \epsilon_2}:X_{i\epsilon_1,s}^{+}({z})X_{j\epsilon_2,s}^{+}({w}):f(z-q^{-1}w)
f(zq^{\frac{\epsilon_1}2}-wq^{\frac{\epsilon_2}2})^{\epsilon_1\epsilon_2}(-1)^{\frac{1-\epsilon_1}2}.
\end{eqnarray*}

Using \eqref{fcn1}-\eqref{fcn2} we have that
\begin{eqnarray*}
&&(z-q^{-\frac12}w)X_{i,s}^{+}({z})X_{j,s}^{+}({w})=:X_{i+}^{+}X_{j+}^{+}:q^{-\frac{1}4}-:X_{i-}^{+}X_{j-}^{+}:q^{\frac{1}4}\\
&&\hskip 50pt +:X_{i+}^{+}X_{j-}^{+}:(q^{\frac14}z-q^{-\frac14}w)-:X_{i-}^{+}X_{j+}^{+}:(q^{-\frac14}z-q^{\frac14}w)\\
&&=(q^{-\frac12}z-w)X_{j,s}^{+}({w})X_{i,s}^{+}({z}),
\end{eqnarray*}
where we have used the fact that $:X_{i\epsilon_1}^{+}(z)X_{j\epsilon_2}^{+}(w):=-:X_{j\epsilon_2}^{+}(w)X_{i\epsilon_1}^{+}(z):$
or $\varepsilon_i\varepsilon_j=-\varepsilon_j\varepsilon_i$.

\subsection{Proof of relation $(2.9)$.}
To prove $(2.9)$ is to check that $$[\,X_{i,s}^+({z}),X_{j,s}^-({w})\,]=\frac{\delta_{ij}}{q_i-q^{-1}_i}\Big(\Phi_i^{(s)}(wq^{\frac{1}2})\delta(\frac{wq}{z})
-\Psi_i^{(s)}(wq^{-\frac{1}{2}})\delta(\frac{wq^{-1}}{z})\Big).$$

Note that the case of $(\alpha_i|\alpha_i)=2$ was known in type $A$, so we will check the cases $(\alpha_i|\alpha_i)=1$ and $(\alpha_i|\alpha_j)=-\frac{1}{2}$ to show the new situation.

For $i\in I_a$, writing $:X^+_{i\epsilon}X^-_{j\epsilon'}:=:X^+_{i\epsilon, s}(z)X^-_{j\epsilon', s}(w):$ etc., by \eqref{3.1}-\eqref{3.2} we have that
\begin{eqnarray*}
&&X_{i,s}^{+}({z})X_{i,s}^{-}({w})=:X_{i+}^{+}X_{i+}^{-}:q^{\frac12}+:X_{i+}^{+}X_{i-}^{-}:(z-w)^{-1}(q^{\frac12}z-q^{-\frac{1}{2}}w)^{-1}\\
&&\hspace{1cm}+:X_{i-}^{+}X_{i+}^{-}:
(z-w)^{-1}(q^{-\frac12}z-q^{\frac{1}{2}}w)^{-1}
+:X_{i-}^{+}X_{i-}^{-}:q^{-\frac{1}{2}}.
\end{eqnarray*}

Modifying the normal order to $\bullet X^{\pm}_{i}(z)X^{\pm}_{j}(w)\bullet$ by excluding $zw$, we get that
\begin{eqnarray*}
[X_{i,s}^{+}({z}), X_{i,s}^{-}({w})]&=&\bullet X_{i+,s}^{+}({z})X_{i-,s}^{-}({w})\bullet\frac{1}{q^{\frac{1}{2}}-q^{-\frac{1}{2}}}(\delta(\frac{w}{z})-\delta(\frac{q^{-1}w}{z}))\\
&&\hspace{0.5cm}+\bullet X_{i-,s}^{+}({z})X_{i+,s}^{-}({w})\bullet\frac{1}{q^{\frac{1}{2}}-q^{-\frac{1}{2}}}(\delta(\frac{qw}{z})-\delta(\frac{w}{z}))\\
&=&\frac{1}{q_i-q_i^{-1}}\bigg(\Phi^{(s)}_i(wq^{\frac{1}{2}})\delta(\frac{qw}{z})-\Psi^{(s)}_i(wq^{-\frac{1}{2}})\delta(\frac{q^{-1}w}{z})\bigg),
\end{eqnarray*}
where we have used the property of the $\delta$-function and the fact that
\begin{eqnarray*}
&&\bullet X_{i+,s}^{+}({z})X_{i-,s}^{-}({qz})\bullet=\Phi_i^{(s)}(q^{1/2}z), \quad
\bullet X_{i-,s}^{+}({z})X_{i+,s}^{-}({q^{-1}z})\bullet=\Psi_i^{(s)}(q^{-1/2}z), \\
&& \bullet X_{i+,s}^{+}({z})X_{i-,s}^{-}({z})\bullet=\bullet X_{i-,s}^{+}({z})X_{i+,s}^{-}({z})\bullet.
\end{eqnarray*}

For the case $(\alpha_i|\alpha_j)=-\frac{1}{2}$, recalling $\eqref{3.2}$ and \eqref{fcn1} we obtain that
\begin{eqnarray*}
X_{i,s}^{+}({z})X_{j,s}^{-}({w})&=&:X_{i+}^{+}X_{j-}^{-}:(q^{\frac14}z-q^{-\frac14}w)-:X_{i+}^{+}X_{j+}^{-}:(q^{-\frac14}z-q^{\frac14}w)\\
&& \hskip 20pt +:X_{i-}^{+}X_{j-}^{-}:q^{-\frac14}
-:X_{i-}^{+}X_{j+}^{-}:q^{\frac{1}{4}}.
\end{eqnarray*}
\begin{eqnarray*}
X_{j,s}^{-}({w})X_{i,s}^{+}({z})&=&:X_{j-}^{-}X_{i+}^{+}:(q^{-\frac14}w-q^{\frac14}z)-:X_{j+}^{-}X_{i+}^{+}:(q^{\frac14}w-q^{-\frac14}z)\\
&& \hskip 20pt +:X_{j-}^{-}X_{i-}^{+}:q^{\frac14}
-:X_{j+}^{-}X_{i-}^{+}:q^{-\frac{1}{4}}.
\end{eqnarray*}

Thus for $(\alpha_i\alpha_j)=-1/2$, we have $[X_{i,s}^{+}({z}), X_{j,s}^{-}({w})]=0$ due to $\varepsilon_i\varepsilon_j=-\varepsilon_j\varepsilon_i$.

\subsection{Proof of relation $(2.10)$}
We will check some cases of $a_{0,1}=-1=a_{n,n-1}$,
as the other cases are similar. By \eqref{3.1}-\eqref{3.2} we have that 
\begin{eqnarray*}
&&X_{0,s}^{+}(z_1)X_{0,s}^{+}(z_2)X_{1\epsilon,s}^{+}(w)\\
&&\hspace{-0.5cm}=:X_{0,s}^{+}(z_1)X_{0,s}^{+}(z_2)X_{1\epsilon,s}^{+}(w):\frac{(z_1-z_2)(z_1-q^{-2}z_2)}{(z_1-q^{-1}w)(z_2-q^{-1}w)}.
\end{eqnarray*}

Subsequently 
\begin{eqnarray*}
&&X_{0,s}^{+}(z_1)X_{0,s}^{+}(z_2)X_{1\epsilon,s}^{+}(w)-(q^{}+q^{-1})X_{0,s}^{+}(z_1)X_{1\epsilon,s}^{+}(w)X_{0,s}^{+}(z_2)+X_{1\epsilon,s}^{+}(w)X_{0,s}^{+}(z_1)X_{0,s}^{+}(z_2)\\
&&=
F\cdot\bigg((w-q^{-1}z_{1})(w-q^{-1}z_{2})+[2](z_{2}-q^{-1}w)(w-q^{-1}z_{1})+(z_{1}-q^{-1}w))(z_{2}-q^{-1}w)\bigg)\\
&&=F\cdot(q^{}-q^{-1})w(z_2-q^{-2}z_1),
\end{eqnarray*}
where $F=\frac{:X_{0,s}^{+}(z_1)X_{0,s}^{+}(z_2)X_{1\epsilon,s}^{+}(w):}{\prod\limits_{i}(z_{i}-q^{-1}w)(w-q^{-1}z_{i})}(z_1-z_2)(z_1-q^{-2}z_2)$.
Observe that $F\cdot (z_2-q^{-2}z_1)$ is antisymmetric in $z_i$, so the Serre relation holds in this case.

When $a_{ij}=-1=a_{ji}\,\,(i, j\in I_a)$, it follows from \eqref{3.1}-\eqref{3.2} that
\begin{eqnarray}\notag
&&\hspace{-0.4cm}X_{i\epsilon_1,s}^{+}(z_1)X_{i\epsilon_2,s}^{+}(z_2)X_{j\epsilon,s}^{+}(w)-[2]_iX_{i\epsilon_1,s}^{+}(z_1)X_{j\epsilon,s}^{+}(w)X_{i\epsilon_2,s}^{+}(z_2)
+X_{j\epsilon,s}^{+}(w)X_{i\epsilon_1,s}^{+}(z_1)X_{i\epsilon_2,s}^{+}(z_2)\\ \notag
&&\hspace{-0.5cm}=:X_{i\epsilon_1,s}^{+}(z_1)X_{i\epsilon_2,s}^{+}(z_2)X_{j\epsilon,s}^{+}(w):
(z_1-q^{-1}z_2)(z_1q^{\epsilon_1/2}-z_2q^{\epsilon_2/2})^{\epsilon_1\epsilon_2}\\ \notag
&&\left(f(z_1-q^{-1}w)f(q^{\frac{\epsilon_1}2}z_1-q^{\frac{\epsilon}2}w)^{\epsilon_1\epsilon}
\cdot f(z_2-q^{-1}w)f(q^{\frac{\epsilon_2}2}z_2-q^{\frac{\epsilon}2}w)^{\epsilon_2\epsilon}
(-1)^{1-\frac{\epsilon_1+\epsilon_2}2} \right.\\ \notag
&&\hskip 0.5cm+[2]_i
f(z_1-q^{-1}w)f(q^{\frac{\epsilon_1}2}z_1-q^{\frac{\epsilon}2}w)^{\epsilon_1\epsilon}
\cdot f(w-q^{-1}z_2)f(q^{\frac{\epsilon}2}w-q^{\frac{\epsilon_2}2}z_2)^{\epsilon_2\epsilon}
(-1)^{1-\frac{\epsilon_1+\epsilon}2}\\ \label{Serre}
&&\hskip 1.5cm\left.
+f(w-q^{-1}z_1)f(q^{\frac{\epsilon_1}2}z_1-q^{\frac{\epsilon}2}w)^{\epsilon_1\epsilon}
\cdot f(w-q^{-1}z_2)f(q^{\frac{\epsilon}2}w-q^{\frac{\epsilon_2}2}z_2)^{\epsilon_2\epsilon}
(-1)^{1-\frac{\epsilon+\epsilon_1}2}\right).
\end{eqnarray}

We claim that the above is antisymmetric with respect to the indices, i.e. under the action: $z_1\to z_2, \epsilon_1\to \epsilon_2$, which
will immediately imply the Serre relation in this case. In fact, there are 8 subcases to verify. We check
one complicated subcase $(\epsilon_1, \epsilon_2, \epsilon)=(+, -, +)$. Then the above expression \eqref{Serre} is simplified to
\begin{eqnarray*}
&&\hspace{-0.3cm}:X_{i+,s}^{+}(z_1)X_{i-,s}^{+}(z_2)X_{j+,s}^{+}(w):\frac{(z_1-q^{-1}z_2)(q^{\frac12}z_1-q^{-\frac12}z_2)^{-1}}
{\prod\limits_{i}(z_i-q^{-\frac12}w)\cdot(w-q^{-\frac12}z_1)}(q^{\frac12}-q^{-\frac12})w(z_2-q^{-1}z_1),
\end{eqnarray*}
where we have used the following identity \cite{JKM}:
\begin{eqnarray*}
&&(w-az_1)(w-az_2)+(a+a^{-1})(z_1-aw)(w-az_2)+(z_1-aw)(z_2-aw)\\
&&\hspace{-0.4cm}=(a^{-1}-a)w(z_1-a^2z_2).
\end{eqnarray*}

Similarly, \eqref{Serre} at $(\epsilon_1, \epsilon_2, \epsilon)=(-, +, +)$ is equal to
\begin{eqnarray*}
&&\hspace{-0.3cm}:X_{i-,s}^{+}(z_1)X_{i+,s}^{+}(z_2)X_{j+,s}^{+}(w):\frac{(z_1-q^{-1}z_2)(q^{-\frac12}z_1-q^{\frac12}z_2)^{-1}}
{\prod\limits_{i}(z_i-q^{-\frac12}w)\cdot(w-q^{-\frac12}z_2)}(q^{\frac12}-q^{-\frac12})w(z_2-q^{-1}z_1).
\end{eqnarray*}
This means that \eqref{Serre} is indeed antisymmetric under $z_1\to z_2, \epsilon_1\to \epsilon_2$. Therefore we have proved that
for $i, j\in I_a, (\alpha_i|\alpha_j)=-1/2$
\begin{eqnarray*}
&&\hspace{-0.4cm}\underset{1\leftrightarrow 2}{Sym}\left(X_{i\epsilon_1,s}^{+}(z_1)X_{i\epsilon_2,s}^{+}(z_2)X_{j\epsilon,s}^{+}(w)-[2]_iX_{i\epsilon_1,s}^{+}(z_1)X_{j\epsilon,s}^{+}(w)X_{i\epsilon_2,s}^{+}(z_2)
\right.\\
&&\hspace{7cm}\left. +X_{j\epsilon,s}^{+}(w)X_{i\epsilon_1,s}^{+}(z_1)X_{i\epsilon_2,s}^{+}(z_2)\right)=0.
\end{eqnarray*}

Now we consider the case of $a_{1,0}=-2=a_{n-1,n}$. By Prop. 3.1 
we have that 
\begin{eqnarray*}
&&\hspace{-0.3cm}X_{1\epsilon_1,s}^{+}(z_1)X_{1\epsilon_2,s}^{+}(z_2)X_{1\epsilon_3,s}^{+}(z_{3})X_{0,s}^{+}(w)\\
&&\hspace{-0.8cm}=:X_{1\epsilon_1,s}^{+}(z_1)X_{1\epsilon_2,s}^{+}(z_2)X_{1\epsilon_3,s}^{+}(z_{3})X_{0,s}^{+}(w):
\frac{\prod\limits_{i<j}(z_i-q^{-1}z_j)(q^{\epsilon_i/2}z_i-q^{\epsilon_j/2}z_j)^{\epsilon_i\epsilon_j}}
{\prod_i(z_i-q^{-1}w)}.
\end{eqnarray*}
Pulling out the common factor we get that
\begin{eqnarray*}
&&\sum\limits_{k=0}^{3}(-1)^k\Big[{3\atop  k}\Big]_{i}X_{1\epsilon_1, s}^{+}({z_1})\cdots X_{1\epsilon_k, s}^{+}(z_k) X_{0,s}^{+}(w)X_{1\epsilon_{k+1}, s}^{+}({z_{k+1}})\cdots
X_{1\epsilon_3, s}^{+}({z_3})\\
&&\hspace{-0.5cm}=F\bigg((z_1-qw)(z_2-qw)(z_{3}-qw)+[3]_{1}(z_1-qw)(z_2-qw)(w-qz_{3})\\
&&\hspace{0.9cm}+[3]_{1}(z_1-qw)(w-qz_2)(w-qz_{3})+(w-qz_1)(w-qz_2)(w-qz_{3})\bigg).\\
&&\hspace{-0.5cm}=F\cdot(q^{-1}-q)\bigg(w^2(z_1-(q+q^2)z_2+q^3z_{3})+w(z_1z_2-(q+q^2)z_1z_{3}+q^3z_2z_{3})\bigg),
\end{eqnarray*}
where $F=:X^+_{1\epsilon_1}(z_1)\cdots X^+_0(w):\frac{\prod\limits_{1\leq i<j\leq 3}(z_{i}-q^{-1}z_{j})(q^{\epsilon_{i}/2}z_{i}-q^{\epsilon_{j}/2}z_{j})^{\epsilon_{i}\epsilon_{j}}}{\prod\limits_{1\leq i\leq 3}(z_{i}-q^{-1}w)(w-q^{-1}z_{i})}q^{-3}$. Now consider the action of ${\mathfrak S}_3$ defined by
$\sigma (z_i, \epsilon_i)=(z_{\sigma(i)}, \epsilon_{\sigma(i)})$. The factor $F$ excluding $\prod\limits_{1\leq i<j\leq 3}(z_{i}-q^{-1}z_{j})$ is antisymmetric under the action of ${\mathfrak S}_3$. Recall the following identity \cite{JKM}:
\begin{eqnarray*}
&&\sum\limits_{\sigma\in\mathfrak S_3}(-1)^{l(\sigma)}(z_1-(q+q^2)z_2+q^3z_{3})(z_1-q^{-1}z_2)(z_1-q^{-1}z_{3})(z_2-q^{-1}z_{3})=0.
\end{eqnarray*}
Replacing $z_i$ by $(z_1z_2z_3)/z_i$ and factoring out an appropriate factor, we also have the following identity
\begin{eqnarray*}
&&\sum\limits_{\sigma\in\mathfrak S_3}(-1)^{l(\sigma)}(z_2z_3-(q+q^2)z_1z_3+q^3z_1z_3)(z_1-q^{-1}z_2)(z_1-q^{-1}z_{3})(z_2-q^{-1}z_{3})=0.
\end{eqnarray*}

Using these two identities, we have proved the following equation:
\begin{eqnarray*}
&&{Sym}\sum_{k=0}^{m=3}(-1)^k\Big[{m\atop  k}\Big]_{1}X_{1\epsilon_1,s}^{\pm}({z_1})\cdots X_{1\epsilon_k,s}^{\pm}(z_k) X_{0,s}^{\pm}({w})X_{1\epsilon_{k+1},s}^{\pm}({z_{k+1}})\cdots X_{1\epsilon_{m},s}^{\pm}({z_m})=0,
\end{eqnarray*}
where $Sym$ runs through the action of $\mathfrak{S}_m$ on $(z_i, \epsilon_i)$, and this then implies $(2.10)$.

For the last relation, without loss of generality, we take $i\in I_a$ in the ``+'' case as an example. Note that for $s\neq s'$
{\begin{eqnarray*}
&&X^+_{i\epsilon_1,s}(z_1)X^+_{i\epsilon_2,s}(z_2)X^+_{i\epsilon_3,s}(z_3)X^-_{i\epsilon,s'}(w)
=:X^+_{i\epsilon_1,s}(z_1)X^+_{i\epsilon_2,s}(z_2)X^+_{i\epsilon_3,s}(z_3)X^-_{i\epsilon,s'}(w):\\
&&\hspace{1.5cm}\cdot\prod\limits_{1\leq i<j\leq 3}(z_i-z_jq^{-1})(z_iq^{\epsilon_i/2}-z_jq^{\epsilon_j/2})^{\epsilon_i\epsilon_j}
(z_1z_2z_3)^{-1}(z_1^{\epsilon_1}z_2^{\epsilon_2}z_3^{\epsilon_3})^{\epsilon}. 
\end{eqnarray*}

When we move $X^-_{i\epsilon,s'}(w)$ across $X^+_{i\epsilon_i}(z_i)$ to the left, the contraction function only changes in the
last three factors respectively to $(z_1z_2w)^{-1}(z_1^{\epsilon_1}z_2^{\epsilon_2}w^{\epsilon_3})^{\epsilon}$, 
$(z_1w^2)^{-1}(z_1^{\epsilon_1}w^{\epsilon_2+\epsilon_3})^{\epsilon}$, 
and 
$(w^3)^{-1}(w^{\epsilon_1+\epsilon_2+\epsilon_3})^{\epsilon}$. 
}

Then we have that
\begin{eqnarray*}
&&\lim_{z_k\to w}\sum\limits_{\epsilon,\epsilon_j\in\{\pm\}}\Bigl(X^+_{i\epsilon_1,s}(z_1)X^+_{i\epsilon_2,s}(z_2)
X^+_{i\epsilon_3,s}(z_3)X^-_{i\epsilon,s'}(w)\\
&&\hspace{3.0cm}-[3]_{i}X^+_{i\epsilon_1,s}(z_1)X^+_{i\epsilon_2,s}(z_2)X^-_{i\epsilon,s'}(w)X^+_{i\epsilon_3,s}(z_3)\\
&&\hspace{4.5cm}+[3]_{i}X^+_{i\epsilon_1,s}(z_1)X^-_{i\epsilon,s'}(w)X^+_{i\epsilon_2,s}(z_2)X^+_{i\epsilon_3,s}(z_3)\\
&&\hspace{6.0cm}-X^-_{i\epsilon,s'}(w)X^+_{i\epsilon_1,s}(z_1)X^+_{i\epsilon_2,s}(z_2)X^+_{i\epsilon_3,s}(z_3)\Bigr)\\
&=&\lim_{z_k\to w}\sum\limits_{\epsilon,\epsilon_j\in\{\pm\}}:X^+_{i\epsilon_1,s}(z_1)X^+_{i\epsilon_2,s}(z_2)X^+_{i\epsilon_3,s}(z_3)
X^-_{i\epsilon,s'}(w):\prod\limits_{1\leq i<j\leq 3}(z_iq^{\epsilon_i/2}-z_jq^{\epsilon_j/2})^{\epsilon_i\epsilon_j}\\
&&\hspace{0.5cm}\times(z_i-z_jq^{-1})
 \Bigl((z_1z_2z_3)^{-1}(z_1^{\epsilon_1}z_2^{\epsilon_2}z_3^{\epsilon_3})^{\epsilon} 
 -[3]_i(z_1z_2w)^{-1}(z_1^{\epsilon_1}z_2^{\epsilon_2}w^{\epsilon_3})^{\epsilon}\\    
&&\hspace{2.0cm}+[3]_i(z_1w^2)^{-1}(z_1^{\epsilon_1}w^{\epsilon_2+\epsilon_3})^{\epsilon}   
 -(w^3)^{-1}(w^{\epsilon_1+\epsilon_2+\epsilon_3})^{\epsilon}\Bigr)=0.  
\end{eqnarray*}}

\vskip30pt \centerline{\bf ACKNOWLEDGMENT}

N. Jing would like to thank the support of
Simons Foundation grant 523868 and NSFC grant 11531004.  H. Zhang would
like to thank the support of NSFC grant 11871325.

\bibliographystyle{amsalpha}

\begin{thebibliography}{99}

\bibitem
{BM} S.~Berman and R.~V.~Moody, \textit{ Lie algebras graded by finite root systems and the intersection matrix algebras of Slowdowy}, Invent. Math. \textbf{108} (1992), 323--347.


\bibitem
{CP1} V.~Chari and A.~Pressley, \textit{A guide to quantum groups}, Cambridge Univ.
Press, Cambridge, 1994.


\bibitem
{CP3} V.~Chari and A.~Pressley, \textit{Quantum affine algebras and their representations}. Representations of groups (Banff, AB, 1994), CMS Conf. Proc., \textbf{16}, Amer. Math. Soc. Providence, RI (1995), 59--78.




\bibitem
{FJM1} B. Feigin, M. Jimbo, T. Miwa and E. Mukhin, \textit{Branching rules for quantum toroidal $\mathfrak{gl}(n)$},
Adv. Math. \textbf{300} (2016), 229-274. 

\bibitem
{FJM2} B. Feigin, M. Jimbo, T. Miwa and E. Mukhin, \textit{Representations of quantum toroidal $\mathfrak{gl}_n$},
J. Algebra \textbf{380} (2013), 78--108. 

\bibitem
{FJ} I.~Frenkel and N.~Jing, \textit{Vertex representations of quantum affine algebras}, Proc. Nat'l. Acad. Sci.
USA. \textbf{85} (1998), 9373--9377.

\bibitem
{FJW} I.~B.~Frenkel, N.~Jing and W.~Wang,
 \textit{Quantum vertex representations via finite groups and the McKay correspondence}, Comm. Math. Phys. \textbf{211} (2000), 365--393.

\bibitem
{FK} I.~Frenkel and V.~G.~Kac, \textit{Basic representations of affine Lie algebras and dual resonance models}, Invent. Math. \textbf{62} (1980), 23-66.

\bibitem
{GJ} Y.~Gao and N. Jing,
\textit{$U_q(\mathfrak{gl}_N)$ action on $\mathfrak{gl}_N$-modules and quantum toroidal algebras}, J. Algebra \textbf{273} (2004), 320--343.

 \bibitem
 {GJXZ}
 Y.~Gao, N.~Jing, L. Xia and H. Zhang, \textit{Quantum $N$-toroidal algebras and quantized GIM algebras of $N$-fold affinization}, 2019, arXiv: 1907.06301.

\bibitem
{GTL} S.~Gautam and V.~Toledano-Laredo,
\textit{Yangians and quantum loop algebras},
Selecta Math. (N.S.) \textbf{19} (2013), 271--336.

\bibitem
{GKV} V.~Ginzburg, M.~Kapranov and E.~Vasserot, \textit{Langlands reciprocty for algebric surfaces,}  Math. Res. Lett. \textbf{2} (1995), 147--160.

\bibitem
{GM} N. Guay and X. Ma, \textit{From quantum loop algebras to Yangians},
J. Lond. Math. Soc. 
\textbf{86} (2012), 683--700.

\bibitem
{H1} D.~Hernandez, \textit{Drinfeld coproduct, quantum fusion tensor category and applications}, Proc. London Math. Soc. \textbf{95}
(2007), 567--608.


\bibitem
{H2} D.~Hernandez, \textit{Quantum toroidal algebras and their representations},
Selecta Math. (N.S.) \textbf{14} (2009), 701--725.

\bibitem
{J1} N.~Jing, \textit{Twisted vertex representations of quantum affine algebras},
Invent. Math. \textbf{102} (1990), 663--690.


\bibitem
{J3} N.~Jing, \textit{ Higher level representations of the quantum affine algebra
$U_q(\hat{sl}(2))$}, J. Algebra {\bf 182} (1996), 448--468.

\bibitem
{J4} N.~Jing, \textit{Quantum Kac-Moody algebras and vertex representations},
Lett. Math. Phys. \textbf{4} (1998), 261--271.

\bibitem
{JKK} N.~Jing, S.-J.~Kang and Y.~Koyama,
\textit{Vertex operators between level one irreducible representations of the
quantum affine algebra $U_q(D^{(1)})$}, Comm. Math. Phys. {\bf 174} (1995),
367--392.

\bibitem
{JKM} N.~Jing, Y.~Koyama and K.~C.~Misra, \textit{Level one representations of quantum affine algebras $U_q(C_n^{(1)})$},
Sel. Math. New Ser.  \textbf{5} (1999), 243--255.


\bibitem
{JM} N.~Jing and K.~C.~Misra, \textit{Vertex operators for twisted quantum affine algebras},
Trans. Amer. Math. Soc. \textbf{351} (1999), 1663--1690.

\bibitem
{Ko} S.~Kolb, \textit{Quantum symmetric Kac-Moody pairs}, Adv. Math. \textbf{267} (2014), 395--469.

\bibitem
{M1} K. Miki, \textit{ Toroidal and level 0 $U_q'(\widehat{sl}_{n+1})$ actions on $U_q(\widehat{gl}_{n+1})$-modules} J. Math. Phys. \textbf{40} (1999), 3191--3210.

\bibitem
{M2} K. Miki, \textit{Representations of quantum toroidal algebra $U_q(sl_{n+1},tor)(n>2)$}, J. Math. Phys. \textbf{41} (2000), 7079--7098.

\bibitem
{M3} K. Miki, \textit{Quantum toroidal algebra $U_q(sl_2,tor)$ and R matrices}, J. Math. Phys. \textbf{42} (2001), 2293--2308.

\bibitem
{MRY}
 R.~V~Moody, S.~E.~Rao and T.~Yokonuma, \textit{Toroidal Lie algebras and vertex representations}, Geom. Ded. \textbf{35} (1990), 283-307.



\bibitem
{RM} S.~Rao and R. Moody, \textit{Vertex representations for $N$-toroidal Lie algebras and a generalization of the Virasoro algebra}, Comm. Math. Phys.  \textbf{159} (1994), 239--264.



\bibitem
{S} Y.~Saito, \textit{Quantum toroidal algebras and their vertex representations}, Publ. RIMS. Kyoto Univ.  \textbf{34} (1998), 155--177.

\bibitem
{Sl} P.~Slodowy, \textit{Beyond Kac-Moody algebras, and inside}. Lie algebras and related topics (Windsor, Ont., 1984), pp. 361--371,
CMS Conf. Proc., \textbf{5}, Amer. Math. Soc., Providence, RI, 1986.


\bibitem
{VV} M.~Varagnolo and E.~Vasserot, \textit{Schur duality in the toroidal setting},  Comm. Math.
Phys. \textbf{182} (1996), 469--484.


\end{thebibliography}

\end{document}